\documentclass[11pt]{article}
\usepackage{epsfig}
\usepackage{graphicx}
\usepackage{amssymb}
\usepackage{amsmath}
\usepackage{bm}
\usepackage{mathrsfs}
\newcommand{\be}{\begin{equation}}
\newcommand{\ee}{\end{equation}}
\newcommand{\bea}{\begin{eqnarray}}
\newcommand{\eea}{\end{eqnarray}}

\usepackage{bm}
\usepackage{color}
\usepackage{graphicx}
\usepackage{cite}

\begin{document}


\title{Unbounded knapsack problem and double partitions}
\author{Boris Y. Rubinstein\\
Stowers Institute for Medical Research
\\1000 50th St., Kansas City, MO 64110, U.S.A.}
\date{\today}

\maketitle
\begin{abstract}
The unbounded knapsack problem can be considered as a particular case
of the double partition problem that asks for a number of nonnegative
integer solutions to a system of two linear Diophantine
equations with integer coefficients.
In the middle of 19th century Sylvester and Cayley
suggested an approach based on the variable elimination allowing a reduction 
of a double partition to a sum of scalar
partitions.
This manuscript discusses a geometric interpretation of this method and 
its application to the knapsack problem.
\end{abstract}

{\bf Keywords}: unbounded knapsack problem, double partition.

{\bf 2010 Mathematics Subject Classification}: 11P82.

\section{Integer partitions and knapsack problem}
\label{intro}

The intimate relation between the knapsack problem and integer partitions 
attracted attention of many researchers \cite{Bradley1971, Faaland1973, Horowitz1974}
who mainly focused on the algorithmic approaches 
while the partition computational aspect was usually neglected. 
It happened probable because in the field 
of partitions there existed no explicit expressions for both scalar and 
vector partitions (to the latter the knapsack problem can be reduced).
The author of this manuscript provided a closed formula for 
evaluation of scalar partitions based on the Bernoulli polynomials 
of higher order \cite{Rub04}. 
In the middle of the 19th century Sylvester suggested an algorithm
of reduction of the vector partition to the sum of scalar ones  \cite{Sylv1}
and Cayley successfully implemented it for the case of double partitions  \cite{Cayley1860}.
This manuscript discusses the application of the original Cayley method to
the unbounded knapsack problem and provides a simple geometric interpretation 
of the algorithm.

\subsection{Scalar partitions}
\label{Scalar}
The problem of integer partition into a set of integers 
is equivalent to a problem of number
of nonnegative integer solutions of the Diophantine equation
\be
s = \sum_{i=1}^m x_i d_i = {\bf x}\cdot{\bf d}.
\label{coin1}
\ee
A scalar
partition function $W(s,{\bf d}) \equiv W(s,\{d_1,d_2,\ldots,d_m\})$
solving the above problem is a
number of partitions of an integer $s$ into positive integers
$\{d_1,d_2,\ldots,d_m\}$. 
The generating function for $W(s,{\bf d})$
has a form
\be
G(t,{\bf d})=\prod_{i=1}^m\frac{1}{1-t^{d_{i}}}
 =\sum_{s=0}^{\infty} W(s,{\bf d})\;t^s\;,
\label{WGF}
\ee
Sylvester proved \cite{Sylv2} a statement about splitting of the scalar partition
 into periodic and non-periodic parts and showed that it may be presented as a sum of 
 the so called Sylvester waves 
\be
W(s,{\bf d}) = \sum_{n=1} W_n(s,{\bf d})\;,
\label{SylvWavesExpand}
\ee
where summation runs over all distinct factors
of the elements of the generator vector ${\bf d}$.
The wave $W_n(s,{\bf d})$ is a quasipolynomial in $s$ 
given by a product of a polynomial multiplied by a 
periodic function with period $n$.
It was shown \cite{Rub04} that it is possible to
express the wave $W_n(s,{\bf d})$ as a finite sum of the Bernoulli polynomials of
higher order.

\subsection{Vector partitions}
\label{Vector}
Consider a function $W({\bf s},{\bf D})$ counting the number of integer
nonnegative
solutions ${\bf x} \ge 0$
to a linear system ${\bf s} = {\bf D} \cdot {\bf x}$, where ${\bf s} = \{s_1,s_2,\ldots,s_l\}$ and
${\bf D}= \{{\bf c}_1,{\bf c}_2,\ldots,{\bf c}_m\}$ is a nonnegative integer $l \times m$ generator matrix
made of $m$ columns  ${\bf c}_i =\{c_{i1},c_{i2},\ldots,c_{il}\},\ (1 \le i \le m \le l)$ where 
some elements $c_{ik}$ might equal zero.
The function $W({\bf s},{\bf D})$ called {\it vector partition} is
a natural generalization of scalar partition to the
vector argument.
The generating function for $W({\bf s},{\bf D})$ reads
\be
G({\bf t},{\bf D})=\prod_{i=1}^m \frac{1}{1-{\bf t}^{{\bf c}_i}} =
\sum_{{\bf s}} W({\bf s},{\bf D}) {\bf t^s},
\quad
{\bf t^s} = \prod_{k=1}^l t_k^{s_k},
\quad
{\bf t}^{{\bf c}_i} = \prod_{k=1}^l t_k^{c_{ik}}.
\label{WvectGF}
 \ee

\subsection{Sylvester-Cayley method of vector partition reduction}
\label{Sylvester-Cayley}

The problem of scalar and vector integer partitions has a long history and J.J. Sylvester 
made a significant contribution to its solution.
In addition to the splitting algorithm for scalar partition \cite{Sylv2}
he suggested \cite{Sylv1} 
to reduce vector partition into a sum of scalar partitions. 
The reduction is an iterative process based on elimination of the 
variable $t_i$ in the generating function (\ref{WvectGF}) which is 
equivalent to the elimination of $x_i$ from the system ${\bf s} = {\bf D} \cdot {\bf x}$.
Sylvester considered a specific double partition problem as 
an illustration of his method and determined regions (chambers) 
each characterized by a unique
expression for vector partition valid in this region only.
He showed that the 
expressions in the adjacent chambers coincide at their
common boundary. 

This approach was successfully applied by 
A. Cayley \cite{Cayley1860} to double partitions subject
to restrictions on the elements of the matrix ${\bf D}$ --
the vectors ${\bf c}_i$ should be noncollinear and
the elements of every column ${\bf c}_i$ must be relatively prime -- which were first 
specified by Sylvester in his lectures \cite{Sylv0}.
The detailed description of the Cayley algorithm derivation is 
presented in \cite{RubDouble2023, RubGen2025} where this method is 
generalized in order to perform the reduction to scalar partitions 
in the cases when the conditions imposed on the  matrix ${\bf D}$ elements
are not satisfied.
In this manuscript we consider the application of the original Cayley algorithm to
the unbounded knapsack problem.



\subsection{Knapsack problem as double partition}
The unbounded knapsack problem (UKP) reads \cite{Martello1990}
\be
\mbox{maximize} \quad V = \sum_{i=1}^m v_i x_i,
\quad
\mbox{subject to} \quad  \sum_{i=1}^m w_i x_i \le W,
\label{UKPorig}
\ee
for nonnegative integer $x_i$. Introduce a slack variable $x_0$ with the corresponding values $w_0=1,v_0=0$ and 
rewrite the problem (\ref{UKPorig}) as a system of the two Diophantine equations
\be
\sum_{i=0}^m w_i x_i = W,
\quad\quad  
\sum_{i=0}^m v_i x_i = V,
\label{UKP2}
\ee
to find $M(W) = \mbox{max}\, V$ for the fixed value of $W$.
We also assume that the nonnegative integers $w_i, v_i$ satisfy the conditions
relation
\be
\gamma_m > \gamma_{m-1} >\ldots \gamma_2 >\gamma_1 >\gamma_0 = 0,
\quad 
\gamma_i = v_i/w_i,
\quad
\mbox{gcd}(v_i,w_i)=1,
\quad
0 \le i \le m,
\label{nocollinear}
\ee 
where $\mbox{gcd}(v_i,w_i)$ denotes the greatest common divisor of $v_i$ and $w_i$.
The system of the two Diophantine equations (\ref{UKP2}) corresponds to the problem 
${\bf s} = {\bf D} \cdot {\bf x}$ with 
${\bf D} = \{{\bf c}_{0},{\bf c}_{1},{\bf c}_{2}, \ldots {\bf c}_{m-1},{\bf c}_{m} \}$ and 
${\bf s} = \{W,V\}^T$, where ${\bf c}_{0} = \{1,0\}^T$, ${\bf c}_{i} = \{w_i,v_i\}^T,\ (1 \le i \le m),$  
are the two-dimensional columns of the matrix ${\bf D}$ and ${}^T$ denotes the transposition.

The UKP can be stated as follows -- for the fixed value of $W$ find the 
maximal value $M(W)$ of $V$ for which the double partition $W({\bf s},{\bf D}) > 0$.
The condition $\mbox{gcd}(w_i,v_i)=1$ means that all columns ${\bf c}_{i}$ are prime ones
while the inequalities in (\ref{nocollinear}) guaranties that the matrix ${\bf D}$ has {\it no collinear columns}.
These conditions coincide with those specified for application of the variable elimination algorithm.
Thus the value of the double partition $W({\bf s},{\bf D})$ can be expressed as a sum of 
scalar partitions that in their turn can be computed through the Bernoulli polynomials of higher order.
This means that the UKP solution can be found as a closed formula 
with the finite number of quasipolynomials terms.



\section{Double partition reduction by variable elimination}
\label{Cayley}

The process of variable elimination in a system of linear equations
is well known and is widely used in the high school algebra course
where it is applied for the solution of a system of $m$ equations
with $m$ unknowns. In this section we describe a generalization
of the elimination method for a system of two Diophantine equations
with $m+1$ variables $x_i, \ 0 \le i \le m$.
Write the system (\ref{UKP2})  as follows
\be
{\bf s} = \sum_{i=0}^m {\bf c}_{i} x_i,
\label{UKP2vect}
\ee
and consider the elimination procedure of the variable $x_j$.
Introduce the integer vector ${\bf n}_{j}$ normal to ${\bf c}_{j}$ satisfying the 
relations (see Fig.\ref{fig1})
\be
{\bf n}_{j} = \{v_j,-w_j\},
\quad
{\bf n}_{j} \cdot {\bf c}_{j}=0,
\quad
{\bf n}_{j} \times {\bf c}_{j}=1,
\label{n_j}
\ee
where ${\bf a} \cdot {\bf b}$ and ${\bf a} \times {\bf b}$ denote the
scalar and vector product of two vectors respectively.
Multiply the vector equation (\ref{UKP2vect}) by ${\bf n}_{j}$ to obtain
\be
l_j = \sum_{i=0}^m d_{ji} x_i,
\quad
l_j =  {\bf s} \cdot{\bf n}_{j} = Wv_j - Vw_j,
\quad
d_{ji} = {\bf c}_{i} \cdot{\bf n}_{j} = w_i v_j - v_i w_j, 
\quad
d_{jj} \equiv 0.
 \label{elim_xj}
\ee
The equation (\ref{elim_xj}) gives rise to the scalar partition
$W_j({\bf s})=W(l_j,{\bf d}_j)$ with the set 
${\bf d}_j=\{d_{j0},d_{j1},d_{j2},\ldots,d_{jm}\}$ of $m$ 
positive integers $d_{ji}, \ i \ne j$.
The theorem proved by Cayley \cite{Cayley1860} states that the value of the 
double partition $W({\bf s},{\bf D})$ is given be the sum of $(m+1)$ scalar partitions
\be
W({\bf s},{\bf D}) = \sum_{j=0}^m W_j({\bf s}) = \sum_{j=0}^m W(l_j,{\bf d}_j).
\label{CayleyRes}
\ee

The elements of the vector ${\bf d}_j$ as well as the integer $l_j$ might be either
positive or negative while  $l_j$ might also vanish.
From the definition of the scalar partition it follows that $W(s,{\bf d})=0$ for the negative $s<0$ and $W(0,{\bf d})=1$.
One has to find a way to evaluate the partition $W(l_j,{\bf d}_j)$ in case when several $d_{ji}$ are negative. 
Assume that $d_{ji}<0$ for $0 \le i \le k < m$,
write the generating function of the scalar partition (\ref{WGF})  
$$
G(t,{\bf d}_j)=\prod_{i=0}^m\frac{1}{1-t^{d_{ji}}},
\quad
i \ne j,
$$
and use the relation
\be
\frac{1}{1-t^{-a}} = - \frac{t^a}{1-t^a},
\label{negd}
\ee
to rewrite $G(t,{\bf d}_j)$ as follows
\be
G(t,{\bf d}_j)=(-1)^{k+1}\prod_{i=0}^k t^{|d_{ji}|} \cdot \prod_{i=0}^m\frac{1}{1-t^{|d_{ji}|}},
\quad
i \ne j.
\label{GF_pos_only}
\ee
Using the definition (\ref{WGF}) of scalar partition as the coefficent of the generating function expansion
we obtain
\be
W(l_j,{\bf d}_j) = (-1)^{k+1} W(l_j + \sum_{i=0}^k d_{ji}, |{\bf d}_j|),
\quad
|{\bf d}_j| = \{|d_{ji}|\}.
\label{scalar_evaluate}
\ee

\section{Sum of scalar partitions for knapsack problem solution}

To apply the result of the double partition reduction into the sum 
of scalar partitions we first have to analyze the individual terms in (\ref{CayleyRes}). 

\subsection{Scalar partition analysis and evaluation}

\begin{figure}[h]
\begin{center}
\includegraphics[width=13.0cm]{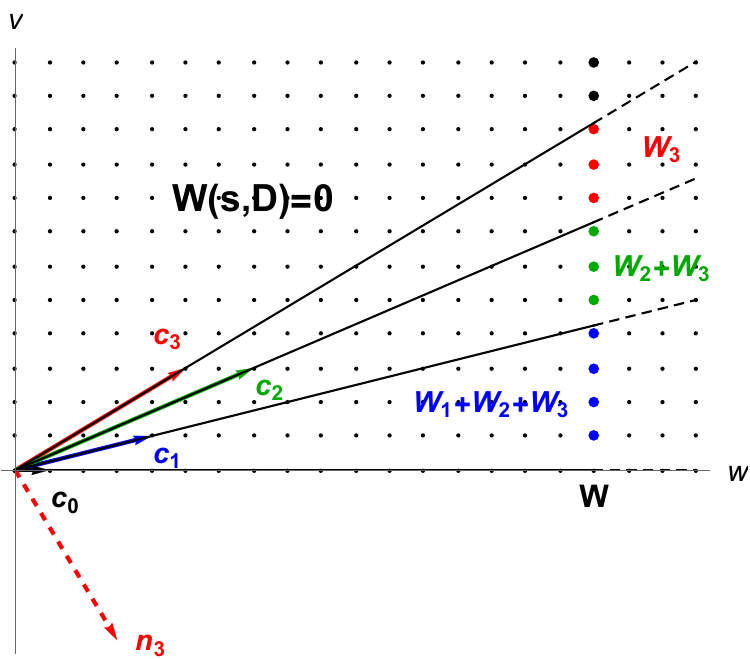}
\end{center}
\caption{
The arrows show the vectors ${\bf c}_i$ for the knapsack problem:  ${\bf c}_0 = \{1,0\}$ (black), ${\bf c}_1 = \{4,1\}$ (blue), 
${\bf c}_2 = \{7,3\}$ (green), ${\bf c}_3 = \{5,3\}$ (red); the normal ${\bf n}_3$ is shown by red dashed arrow.
The double partition $W({\bf s},{\bf D})$ vanishes at all points
that lie above the ray extending the vector ${\bf c}_3$. 
To determine $M(W)$ one has to evaluate  $W({\bf s},{\bf D})$ only at the integer points
${\bf s} = \{W,V\}$  with $0 < V \le \lfloor W \gamma_3 \rfloor$.
The nonzero contributions of $W_j({\bf s})$ into $W({\bf s},{\bf D})$  are shown
in the corresponding chambers. The value of $W({\bf s},{\bf D})$ at the red points
is determined by $W_3({\bf s})$ only, while for the green and blue points one has to 
use the sum $W_2({\bf s})+W_3({\bf s})$ and $W_1({\bf s})+W_2({\bf s})+W_3({\bf s})$, 
respectively. 
}
\label{fig1}
\end{figure}

Considering the element $d_{ji}={\bf c}_{i} \cdot{\bf n}_{j}$ we observe that
this element is positive when the angle between the vectors ${\bf c}_{i}$ and ${\bf n}_{j}$ is less than $\pi/2$ and is negative
when the angle exceeds $\pi/2$. Using the inequalities (\ref{nocollinear}) we find that 
$d_{ji} > 0$ for $j>i$ and $d_{ji} < 0$ for $j<i$ (see Fig.\ref{fig1} for illustration and \cite{Sylv0}). 
Similarly, the sign of $l_j$ is not fixed as it depends
on the angle between ${\bf s}$ and ${\bf n}_{j}$ -- if $ \gamma_j > \Gamma = V/W$ we obtain $l_j > 0$, when 
$\gamma_j < \Gamma$ the value of $l_j$ is negative and finally $l_j=0$ for $\gamma_j = \Gamma$.
Using Fig.\ref{fig1} it is easy to observe that for {\it any} integer point ${\bf s}=\{w,v\}^T$ with $v/w > \gamma_m$ {\it all}
summands $W_j({\bf s})$ equal zero and the double partition $W({\bf s},{\bf D})$ vanishes.
This means that in order to find $M(W)$ one has to test only a finite 
number of points ${\bf s} = \{W,V\}^T$ with $0 < V \le v_{m}$,
where $v_{j}=\lfloor W \gamma_j \rfloor$  and $\lfloor x \rfloor$ denotes integer part of $x$.

Applying this analysis we observe that the term $W_m({\bf s})$ has 
all $d_{mi}$ positive while $l_m > 0$ for $0 < V \le v_{m}$.
This means that the contribution of this term to the double partition is positive.

The term $W_{m-1}({\bf s})$ has a single negative element $d_{m-1,m} < 0$ and 
$l_{m-1} > 0$ for $0 < V \le v_{m-1}$ while 
$l_{m-1} < 0$ for $v_{m-1} < V \le v_{m}$ and in this range 
the single term $W_m({\bf s})$ contributes to the double partition.
We find that the term $W_{m-1}({\bf s})$ contributes negatively
$$
W_{m-1}({\bf s}) = - W(l_{m-1} + d_{m-1,m},  |{\bf d}_{m-1}|).
$$
Similarly we find that the term $W_{m-2}({\bf s})$ generates positive contribution
$$
W_{m-2}({\bf s}) = W(l_{m-2} + d_{m-1,m} + d_{m-2,m-1},  |{\bf d}_{m-2}|).
$$
and it is nonzero for  $0 < V \le v_{m-2}$, so that 
in the range $v_{m-2} < V \le v_{m-1}$
the two terms $W_m({\bf s})$ and $W_{m-1}({\bf s})$ contribute to the double partition with opposite signs.

Continuing this process we find that the last term $W_{0}({\bf s})$ of the expansion (\ref{CayleyRes}) 
has $l_0 < 0$ and all $d_{0i} < 0$ so that it does not contribute into the double partition.
Thus the number of nonzero terms $W_{j}({\bf s})$ equals the number $m$ of the terms in the original problem
(\ref{UKPorig}) and the contribution of these terms is sign alternating.

\subsection{Search algorithm for $M(W)$}

The analysis done in the preceding Sections allows to 
split the interval $I=[1, v_{m}]$ inside which the search should be performed into 
$m$ subintervals $I_{j}=[v_{j-1}+1,v_{j}]$. 
The expression for the double partition $W({\bf s},{\bf D})$ inside the subinterval $I_{j}$
is given by the partial sum
\be
W_j({\bf s},{\bf D}) = \sum_{k=j}^m W_{k}({\bf s}) =
\sum_{k=j}^m (-1)^{m-k} W(l_{j} + \sum_{n=j}^{m-1} d_{n,n+1}, |{\bf d}_{j}| ),
\label{Wpartial}
\ee
with the alternating signs of the consecutive contributions.

As we are searching the maximal value $V$ for which 
$W({\bf s},{\bf D}) \ne 0$ we scan the subintervals 
in the descending order of $V$. 
Note that inside the subinterval all the parameters in (\ref{Wpartial})
are fixed except the value of $l_{j} = Wv_j - V w_j$ that increases with 
the decreasing value of $V$.

We start from the top subinterval $I_{m}$ where $W_m({\bf s},{\bf D})$
is given by a single positive scalar partition $W_{m}({\bf s})$ scanning it starting with 
$V=v_{m}$. The first nonzero value of $W_{m}({\bf s})$ gives $M(W)$.
If the search inside $I_{m}$ fails we move onto $I_{m-1}$ where 
$W_{m-1}({\bf s},{\bf D})$ is provided by the sum of two opposite sign terms.
Thus for each $V$ we evaluate $W_{m}({\bf s})$ and if it is zero we move to the next value,
otherwise we have to find the value of $W_{m-1}({\bf s})$ to check whether it cancels or not
the first term in $W_{m-1}({\bf s},{\bf D})$.
In the subinterval $I_{m-2}$ we first find the sum $W_{m}({\bf s})+W_{m-1}({\bf s})$
and only when it is zero we move on to compute $W_{m-2}({\bf s})$. The procedure repeated
until the value $M(w)$ is determined.


\end{document}